\documentclass[11pt]{article}
\usepackage[english]{babel}
\usepackage{amsmath,amsthm,amssymb}
\usepackage{mathrsfs}
\usepackage{bbm}
\usepackage{amsfonts}
\usepackage[margin=0.8in]{geometry}
\usepackage{graphicx}
\usepackage{caption}
\usepackage{subcaption}
\usepackage{wasysym}
\usepackage{enumerate}
\usepackage{algorithm2e}
\usepackage{tabularx}
\usepackage{color}
\usepackage{wrapfig}
\usepackage{todonotes}

\newtheorem{thm}{Theorem}[section]
\newtheorem{ass}[thm]{Assumption}
\newtheorem{cor}[thm]{Corollary}
\newtheorem{lem}[thm]{Lemma}
\newtheorem{prop}[thm]{Proposition}
\newtheorem*{hyp*}{Hypothesis}
\theoremstyle{definition}
\newtheorem{defn}[thm]{Definition}

\theoremstyle{rem}
\newtheorem{rem}[thm]{Remark}
\numberwithin{equation}{section}

\newcommand{\R}{\mathbb R}

\newcommand{\bbD}{\mathbb D}
\newcommand{\bbF}{\mathbb F}

\newcommand{\bbQ}{\mathbb Q}

\newcommand{\mcA}{\mathcal{A}}

\newcommand{\mcB}{\mathcal{B}}

\newcommand{\mcE}{\mathcal E}

\newcommand{\mcF}{\mathcal F}

\newcommand{\mcT}{\mathcal T}
\newcommand{\mcP}{\mathcal P}

\newcommand{\mcH}{\mathcal H}
\newcommand{\mcK}{\mathcal K}

\newcommand{\mcS}{\mathcal S}
\newcommand{\mcZ}{\mathcal Z}

\newcommand{\E}{\mathbb{E}}
\newcommand{\Prob}{\mathbb{P}}

\newcommand{\esssup}{\mathop{\rm{ess}\,\sup}}

\newcommand{\ett}{\mathbbm{1}}

\newcommand{\cadlag}{c\`adl\`ag~}

\newcommand{\PrM}{\mathfrak{P}}

\newcommand{\ie}{\textit{i.e.\ }}
\newcommand{\eg}{\textit{e.g.\ }}
\newcommand{\etal}{\textit{et.~al.\ }}

\begin{document}


\title{A Note on Reflected BSDEs in Infinite Horizon with Stochastic Lipschitz Coefficients}

\author{Magnus Perninge\footnote{M.\ Perninge is with the Department of Physics and Electrical Engineering, Linnaeus University, V\"axj\"o,
Sweden. e-mail: magnus.perninge@lnu.se.}} %
\maketitle
\begin{abstract}
We consider an infinite horizon, obliquely reflected backward stochastic differential equation (RBSDE). The main contribution of the present work is that we generalize previous results on infinite horizon RBSDEs to the setting where the driver has a stochastic Lipschitz coefficient. As an application we consider robust optimal stopping problems for functional stochastic differential equations (FSDEs) where the driver has linear growth.
\end{abstract}

\section{Introduction}
Backward stochastic differential equations (BSDEs) has been a topic of rapid development during the last decades. Non-linear BSDEs were independently introduced in \cite{PardouxPeng90} and \cite{DuffieEpstein92} and has since found numerous applications.

El Karoui \etal introduced the notion of reflected backward stochastic differential equations (RBSDEs) and demonstrated a link between RBSDEs and optimal stopping in~\cite{ElKaroui1}. This was later extended to infinite horizon in~\cite{HamLepIH} and to discontinuous barriers in \cite{HamRefBSDE}.

Much of the development within the topic of BSDEs have been aimed at showing existence and uniqueness of solutions to BSDEs (reflected and non-reflected) under various types of weaker conditions on the coefficients than those assumed in the original publications. Important contributions from the perspective of the present work are \cite{Bender00,Briand08} that consider BSDEs where the Lipschitz coefficient on the $z$-variable of the driver is a stochastic process and the more recent work presented in~\cite{ElAsri2020} where a reflected BSDE with stochastic Lipschitz coefficient is solved. The objective of the present article is to extend the result and to some extent also the methodology of the latter to an infinite horizon setting.

The method applied to obtain existence of solutions to the infinite horizon RBSDEs considered in the novel paper~\cite{HamLepIH} relies on the approach to infinite horizon (non-reflected) BSDEs developed in \cite{Chen98} and assumes that the Lipschitz coefficients of the driver are deterministic and satisfy an integrability-condition. In particular, this implies that the dependence on $y$ and $z$ in the driver has to vanish as time tends to infinity. We aim to extend this result by establishing existence and uniqueness of solutions to the RBSDE
\begin{align}\label{ekv:bsde-reflection-IH}
  \left\{\begin{array}{l}Y_t=Y_T+\int_t^T f(s,Y_s,Z_s)ds-\int_t^T Z_sdW_s+K_T-K_t,\quad \forall t\in [0,T],\forall T\geq 0,\\
  Y_T\to 0,\: \Prob\mbox{-a.s.~as }T\to\infty\\
  Y_t\geq S_t,\,\forall t\in[0,\infty) \quad\mbox{and}\quad  \int_0^\infty(Y_t-S_t)dK_t=0\end{array}\right.
\end{align}
under the assumption that the driver $f$ satisfies a Lipschitz condition specified in terms of stochastic processes. Our convention of defining infinite horizon bsdes in truncated form is inspired by the work on discrete time reflected backward stochastic difference equations (RBS$\Delta$Es) in~\cite{LondonPaper1} and allows us to consider situations where individual terms on the right hand side diverge as $T\to\infty$.

Notable is also that, as opposed to the setting in for example~\cite{HamLepIH}, the terminal value is set to zero. Although this may appear restrictive it is a natural assumption for many applications, such as problems in optimal control that rarely contain a reward at infinity as a any such reward (however large) gives little comfort at a finite time.

Moreover, we relate solutions to the above reflected BSDE to weak formulations of robust, non-Markovian optimal stopping problems where the aim is to maximize the reward functional
\begin{align}
J(\tau,\alpha):=\E\Big[\int_0^\tau e^{-\rho(t)}\phi(t,(X^{\alpha}_s)_{s\leq t},\alpha_t)dt+e^{-\rho(\tau)}\psi(\tau,(X^{\alpha}_s)_{s\leq \tau})\Big]
\end{align}
over stopping times $\tau$ when simultaneously a minimization if performed over continuous controls $\alpha:=(\alpha_s)_{s\geq 0}$, taking values in a compact subset $A$ of $\R^d$ and $X^{\alpha}$ solves the functional SDE
\begin{align}
X^{\alpha}_t&=x_0+\int_0^ta(s,(X^{\alpha}_r)_{r\leq s},\alpha_s)ds+\int_0^t\sigma(s,(X^{\alpha}_r)_{r\leq s})dW_s.\label{ekv:forward-sde}
\end{align}
In particular, using our more general form of RBSDE \eqref{ekv:bsde-reflection-IH} we are able to relax the common assumption (see \eg \cite{HamLepIH}) that $|\sigma^{-1}(t,x)a(t,x,\alpha)|$ is bounded from above by a square integrable deterministic function and instead assume a linear growth, \ie that $|\sigma^{-1}(t,x)a(t,x,\alpha)|\leq k_L(1+\sup_{s\in[0,t]}|x_s|)$, for some constant $k_L>0$.

In finite horizon, the cooperative version of the above control problem has received a reasonable amount of attention in the last decades. In \cite{Zamfirescu06,Bayraktar2011,ElAsri2020} the problem is solved under various types of assumptions on the involved coefficients. The finite horizon version of the above stochastic differential game of control and stopping (in their terminology) was considered in \cite{Zamfirescu08}, under the assumption that the running and terminal reward are bounded.

The remainder of the article is organized as follows. In the next section we set the notation and recall some well known results for RBSDEs with deterministic Lipschitz coefficients from the original work~\cite{ElKaroui1}. Then in Section~\ref{sec:rbsde-IH} we give as set of assumptions under which we show that \eqref{ekv:bsde-reflection-IH} admits a unique solution. In Section~\ref{sec:robust-impulse} we provide a solution to the above robust optimal stopping problem by relating solutions to \eqref{ekv:bsde-reflection-IH} to weak formulations of the control problem at hand, showing that a saddle-point exists.

\section{Preliminaries\label{sec:prelim}}
We let $(\Omega,\mcF,\bbF,\Prob)$ be a complete filtered probability space, where $\bbF:=(\mcF_t)_{t\geq 0}$ is the augmented natural filtration of a $d$-dimensional Brownian motion $W$ and $\mcF:=\mcF_\infty=\bigcup_{t\geq 0}\mcF_t$.\\

\noindent Throughout, we will use the following notation:
\begin{itemize}
  \item We let $\E$ denote expectation with respect to $\Prob$ and for any other probability measure $\bbQ$ on $(\Omega,\mcF)$, we denote by $\E^{\bbQ}$ expectation with respect to $\bbQ$.
  \item $\mcP_{\bbF}$ is the $\sigma$-algebra of $\bbF$-progressively measurable subsets of $[0,\infty)\times \Omega$.
  \item For $p\geq 1$, we let $\mcS^{p}$ be the set of all $\R$-valued, $\mcP_{\bbF}$-measurable, continuous processes $(Z_t: t\geq 0)$ such that $\|Z\|_{\mcS^p}:=\E\big[\sup_{t\in[0,\infty)} |Z_t|^p\big]^{1/p}<\infty$.
  \item For $p\geq 1$, we let $\mcS^{p,loc}$ be the set of processes that are locally in $\mcS^p$ so that $Z\in \mcS^{p,loc}$ if $Z$ is continuous, $\R$-valued, $\mcP_{\bbF}$-measurable and there is a non-decreasing sequence of $\bbF$-stopping times $(\eta_l)_{l\geq 0}$ with $\eta_l\to\infty$, $\Prob$-a.s., as $l\to\infty$, such that $\|\ett_{[0,\eta_l]}Z\|_{\mcS^p}<\infty$ for all $l\geq 0$.
  \item We let $\mcH^{p}$ denote the set of all $\R^d$-valued $\mcP_{\bbF}$-measurable processes $(Z_t: t\geq 0)$ such that $\|Z\|_{\mcH^p}:=\E\big[\big(\int_0^\infty |Z_t|^2 dt\big)^{p/2}\big]^{1/p}<\infty$. Moreover, we let $\mcH^{p,loc}$ be the local version of $\mcH^p$, so that $Z\in \mcH^{p,loc}$ if there is a non-decreasing sequence of $\bbF$-stopping times $(\eta_l)_{l\geq 0}$ with $\eta_l\to\infty$, $\Prob$-a.s., as $l\to\infty$, such that $\ett_{[0,\eta_l]}Z\in\mcH^p$ for all $l\geq 0$.
  \item We let $\mcT$ be the set of all $\bbF$-stopping times and for each $\eta\in\mcT$ we let $\mcT_\eta$ be the corresponding subsets of stopping times $\tau$ such that $\tau\geq \eta$, $\Prob$-a.s.
  \item We let $\mcA$ be the set of all $\mcP_\bbF$-measurable processes $(\alpha_t)_{t\geq 0}$ taking values in $A$ and for each $t\geq 0$ we let $\mcA_t$ be the set of all $\mcP_\bbF$-measurable processes $(\alpha_s)_{s\geq t}$ taking values in $A$.
  \item We let $*$ denote stochastic integration and set $(X*W)_{t,T}=\int_t^TX_sdW_s$.
  \item We let $\mcE$ denote the Dol\'eans-Dade exponential and use the notation
  \begin{align*}
    \mcE(X*W)_{t,T}=e^{\int_t^TX_sdW_s-\frac{1}{2}\int_t^T|X_s|^2ds}.
  \end{align*}
  Also, we write $\mcE(X*W)_{T}:=\mcE(X*W)_{0,T}$.
  \item For any $T>0$ and any $\mcP_\bbF$-measurable \cadlag process $\zeta$ for which $\E[\mcE(\zeta*W)_T]=1$, we define $\bbQ^{\zeta}$ to be the probability measure equivalent to $\Prob$, such that $d\bbQ^\zeta=\mcE(\zeta*W)_Td\Prob$.
  \item For any non-negative, $\mcP_\bbF$-measurable \cadlag process $L$ we let $\mcZ^L$ be the set of all $\mcP_\bbF$-measurable \cadlag process $\zeta$ with $|\zeta_t|\leq L_t$, for all $t\geq 0$.
  \item For any non-negative, $\mcP_\bbF$-measurable \cadlag process $L$ we let $\PrM^L$ denote the set of all probability measures $\bbQ$ on $(\Omega,\mcF)$ such that $d\bbQ=\mcE(\zeta*W)_Td\Prob$, for some $\zeta\in\mcZ^L$ and some $T>0$.
  \item For any probability measure $\bbQ$ equivalent to $\Prob$, we let $\mcS^p_\bbQ$ and $\mcH^p_\bbQ$ be defined as $\mcS^p$ and $\mcH^p$, respectively, with the exception that the norm is defined with expectation taken with respect to $\bbQ$, \ie $\|Z\|_{\mcS^p_\bbQ}:=\E^\bbQ\big[\sup_{t\in[0,\infty)} |Z_t|^p\big]^{1/p}$ and $\|Z\|_{\mcH^p_\bbQ}:=\E^\bbQ\Big[\big(\int_0^\infty |Z_t|^2 dt\big)^{p/2}\Big]^{1/p}$.
  \item We let $\mcK^p_\bbQ$ be the subset of $\mcS^p_\bbQ$ of non-negative processes.
  \item For any non-negative, $\mcP_\bbF$-measurable \cadlag processes $L^1$ and $L^2$ we let $\mcS^{p,L^1,L^2}$ and $\mcK^{p,L^1,L^2}$ be the subsets of $\cap_{\bbQ\in\PrM^{L^2}}\mcS^{p}_\bbQ$ and $\cap_{\bbQ\in\PrM^{L^2}}\mcK^{p}_\bbQ$, respectively, of all processes $Z$ such that
      \begin{align*}
        \sup_{\bbQ\in \PrM^{L^2}}\|e^{\int_0^\cdot L^1
        _t dt}Z\|_{\mcS^p_\bbQ}<\infty
      \end{align*}
      and
      \begin{align*}
        \sup_{\bbQ\in \PrM^{L^2}}\|\ett_{[T,\infty)}e^{\int_0^\cdot L^1_t dt}Z\|_{\mcS^p_\bbQ}\to 0,
      \end{align*}
       as $T\to\infty$.
  \item For any $\tau\in\mcT$, we add $([0,\tau])$ to the definition of the above spaces to indicate that the space is restricted to processes with index set $[0,\tau]$. For example, $\mcS^p([0,\tau])$ will denote the set of all continuous processes, $(Z_t)_{0\leq t\leq\tau}$, such that $Z_{\cdot\wedge\tau}$ is $\mcP_\bbF$-measurable and $\|Z\|_{\mcS^p([0,\tau])}:=\E\big[\sup_{t\in [0,\tau]} |Z_t|^p\big]^{1/p}<\infty$.
\end{itemize}

In addition, we will throughout assume that, unless otherwise specified, all inequalities hold in the $\Prob$-a.s.~sense.


\subsection{Prior results on RBSDEs}

As our approach will rely heavily on the original work in \cite{ElKaroui1}, we recall the following results with finite horizon $T>0$:

\begin{thm}\label{thm:ElKaroui-rbsde}\emph{(El Karoui \etal \cite{ElKaroui1})}
  Assume that
  \begin{enumerate}[a)]
  \item $\xi\in L^2(\Omega,\mcF_T,\Prob)$
  \item The barrier $S$ is real-valued, $\mcP_\bbF$-measurable and continuous with $S^+\in\mcS^2([0,T])$ and $S_T\leq\xi$.
  \item $f:[0,T]\times \Omega\times \R\times\R^d\to\R$ is such that\footnote{Throughout, we generally suppress dependence on $\omega$ and refer to, for example, $f$ as a map $(t,y,z)\to f(t,y,z)$.} $f(\cdot,y,z)\in\mcH^2([0,T])$ for all $(y,z)\in \R\times\R^{d}$ and for some $k_f>0$ and all $(y,y',z,z')\in\R^{2(1+d)}$ we have
  \begin{align*}
     |f(t,y',z')-f(t,y,z)|\leq k_f(|y'-y|+|z'-z|).
  \end{align*}
  \end{enumerate}
  Then, there exists a unique triple $(Y,Z,K):=(Y_t,Z_t,K_t)_{0\leq t\leq T}$ with $Y,K\in\mcS^2([0,T])$ and $Z\in\mcH^2([0,T])$, where $K$ is non-decreasing with $K_0=0$, such that
\begin{align*}
  \begin{cases}
    Y_t=\xi+\int_t^T f(s,Y_s,Z_s)ds-\int_t^T Z_s dW_s+K_T-K_t,\quad\forall t\in[0,T],\\
    Y_t\geq S_t,\, \forall t\in [0,T] \quad\mbox{and}\quad \int_0^T \left(Y_t-S_t\right)dK_t=0.
  \end{cases}
\end{align*}
Furthermore\footnote{Throughout, $C$ will denote a generic positive constant that may change value from line to line.},
\begin{align}\label{ekv:ElKaroui-bound}
\|Y\|_{\mcS^2([0,T])}^2+\|Z\|_{\mcH^2([0,T])}^2+\|K\|_{\mcS^2([0,T])}^2&\leq C\E\Big[|\xi|^{2}+\int_0^T|f(s,0,0)|^2ds+\sup_{t\in[0,T]}|(S_t)^+|^2\Big].
\end{align}
In addition, $Y$ can be interpreted as the Snell envelope in the following way
\begin{equation*}
      Y_t=\esssup_{\tau\in\mcT_t}\E\bigg[\int_t^\tau f(s,Y_s,Z_s)ds+S_\tau\ett_{[\tau<T]}+\xi \ett_{[\tau=T]}\Big|\mcF_t\bigg]
\end{equation*}
and with $D_t:=\inf\{r\geq t: Y_r=S_r\}\wedge T$ we have the representation
\begin{equation*}
      Y_t=\E\bigg[\int_t^{D_t} f(s,Y_s,Z_s)ds+S_{D_t}\ett_{[D_t<T]}+\xi \ett_{[D_t=T]}\Big|\mcF_t\bigg]
\end{equation*}
and $K_{D_t}-K_t=0$, $\Prob$-a.s.

Moreover, if $(\tilde Y,\tilde Z,\tilde K)$ is the solution to the reflected BSDE with parameters $(\tilde\xi,\tilde f,\tilde S)$, then
\begin{align}\nonumber
&\|\tilde Y-Y\|_{\mcS^2([0,T])}^2+\|\tilde Z-Z\|_{\mcH^2([0,T])}^2+\|\tilde K-K\|_{\mcS^2([0,T])}^2
\\
&\leq C(\|\tilde S-S\|_{\mcS^{2}([0,T])}\Psi_T^{1/2}+\E\Big[|\tilde \xi-\xi|^{2}+\int_0^T |\tilde f(s,Y_s,Z_s)-f(s,Y_s,Z_s)|^2ds\Big]),\label{ekv:ElKaroui-diff}
\end{align}
where
\begin{align*}
\Psi_T:=\E\Big[|\tilde \xi|^{2}+|\xi|^{2}+\int_0^T (|\tilde f(s,0,0)|^2+|f(s,0,0)|^2)ds+\sup_{t\in[0,T]}|(\tilde S_{t})^+ + (S_{t})^+|^{2}\Big].
\end{align*}
\end{thm}


\section{Reflected BSDEs with stochastic Lipschitz coefficient in infinite horizon\label{sec:rbsde-IH}}
We now return to the reflected BSDE from the introduction, namely
\begin{align*}
  \left\{\begin{array}{l}Y_t=Y_T+\int_t^T f(s,Y_s,Z_s)ds-\int_t^T Z_sdW_s+K_T-K_t,\quad \forall t\in [0,T],\forall T\geq 0,\\
  Y_T\to 0,\: \Prob\mbox{-a.s.~as }T\to\infty\\
  Y_t\geq S_t,\,\forall t\in[0,\infty) \quad\mbox{and}\quad  \int_0^\infty(Y_t-S_t)dK_t=0.\end{array}\right.
\end{align*}

Throughout this section we will assume that:
\begin{ass}\label{ass:on-rbsde-IH}
 There are $\mcP_\bbF$-measurable, continuous processes $\bar L^y$, $\underline L^y$ and $L^z$, with $\bar L^y_t\leq \underline L^y_t$, $\Prob$-a.s.~for all $t\in [0,\infty)$, and constants $K_f,K_S>0$, such that:
\begin{enumerate}[(i)]
  \item $\E\big[\mcE(\zeta*W)_T\big]=1$ for all $\zeta\in\mcZ^{L^z}$ and $T\geq 0$.
  \item\label{ass:on-f-IH} The map $(t,\omega,y,z)\mapsto f(t,\omega,y,z): [0,\infty)\times\Omega\times\R\times \R^{d}\to\R$ is such that $(f(t,y,z))_{0\leq t\leq T}\in\mcH^2_\bbQ$ for all $(T,y,z)\in [0,\infty)\times\R\times\R^{d}$ and $\bbQ\in\PrM^{L^z}$. Furthermore, there is a $\bar K^{f}\in \mcK^{2,\bar L^y,L^z}$, with $\|\bar K^{f}\|_{\mcS^2_\bbQ}\leq K_f$ for all $\bbQ\in\PrM^{L^z}$, such that
  \begin{enumerate}[a)]
  \item For all $\bbQ\in\PrM^{L^z}$ and all $t\in[0,\infty)$ we have
  \begin{align}
  \E^\bbQ\Big[\int_{t}^{\infty}e^{\int_t^s\bar L^y_rdr}|f(s,0,0)|ds\Big|\mcF_t\Big]\leq \bar K^{f}_t,
  \end{align}
  $\Prob$-a.s.
  \item The Lipschitz condition
  \begin{align}
     |f(t,y',z')-f(t,y,z)|\leq \underline L^y_t|y'-y|+L^z_t|z'-z|
  \end{align}
  and the growth condition
  \begin{align}
     (f(t,y',z)-f(t,y,z))(y'-y)\leq \bar L^y_t|y'-y|^2
  \end{align}
  hold for all $(t,y,y',z,z')\in[0,\infty)\times\R^{2(1+d)}$, $\Prob$-a.s.
  \end{enumerate}
  \item The barrier $S$ is real-valued, $\mcP_\bbF$-measurable and continuous with $S^+\in\mcS^2_\bbQ$ for all $\bbQ\in\PrM^{L^z}$ and $\limsup_{t\to\infty}S_t\leq 0$, $\Prob$-a.s. Moreover, there is a $\bar K^{S}\in\mcK^{2,\bar L^y,{L^z}}$, with $\|\bar K^{S}\|_{\mcS^2_\bbQ}\leq K_S$ for all $\bbQ\in\PrM^{L^z}$, such that
  \begin{align}
  \E^\bbQ\Big[\sup_{r\in[t,\infty)}e^{\int_t^r\bar L^y_sds}S_r^+\Big|\mcF_t\Big]\leq \bar K^{S}_t,
  \end{align}
  $\Prob$-a.s.~for all $t\in[0,\infty)$ and $\bbQ\in \PrM^{L^z}$.
\end{enumerate}
\end{ass}

The results in this section will rely on applying a different truncation to \eqref{ekv:bsde-reflection-IH} and we introduce the following sequence of stopping times:
\begin{defn}
We let $(\eta_l)_{l\geq 0}$ be the sequence of $\bbF$-stopping times defined as $\eta_l:=\inf\{s\geq 0:\underline L^y_s\vee L^z_s\geq l\}\wedge l$.
\end{defn}
We note that, by continuity of $\underline L^y$ and $L^z$, we have that $\eta_l\to\infty$, $\Prob$-a.s.~as $l\to\infty$. Moreover, on $[0,\eta_l]$ the reflected BSDE \eqref{ekv:bsde-reflection-IH} is of standard Lipschitz type with coefficient $l$.\\

Existence of a solution to \eqref{ekv:bsde-reflection-IH} will be obtained through an approximation scheme and for each $m,n\geq 0$, we introduce the following reflected BSDE,
\begin{align}\label{ekv:trunk-time}
\begin{cases}
  Y^{m,n}_t=\int_{t}^\infty f^{m,n}(s,Y^{m,n}_s,Z^{m,n}_s)ds-\int_t^\infty Z^{m,n}_sdW_s+K^{m,n}_\infty-K^{m,n}_t,\quad\forall t\in[0,\infty),\\
  Y^{m,n}_t\geq S_t,\,\forall t\in[0,\infty)\quad\mbox{and}\quad \int_{0}^\infty(Y^{m,n}_s- S_s)dK^{m,n}_s=0,
\end{cases}
\end{align}
where $f^{m,n}(t,y,z):=\ett_{[0, \eta_m]}(t)f^+(t,y,z)-\ett_{[0,\eta_n]}(t)f^-(t,y,z)$.

\begin{lem}\label{lem:rbsde-trunk-bound-IH}
For each $m,n\geq 0$, the reflected BSDE \eqref{ekv:trunk-time} has a unique solution $(Y^{m,n},Z^{m,n},K^{m,n})\in\mcS^{2}\times\mcH^2\times\mcS^2$, where $K^{m,n}$ is non-decreasing with $K^{m,n}_0=0$. Furthermore, there is a constant $C>0$ (that does not depend on $m,n$ and $\bbQ$) such that
\begin{align}\label{ekv:unif-bound-IH}
\|e^{\int_0^\cdot \bar L^y_rdr}Y^{m,n}\|_{\mcS^2_\bbQ}^2\leq C
\end{align}
and for each $l\geq 0$, there is a $C>0$ (that does not depend on $m,n$ and $\bbQ$, but may depend on $l$) such that
\begin{align}\label{ekv:unif-bound-Z-IH}
\|Z^{m,n}\ett_{[0,\eta_l]}\|_{\mcH^2_\bbQ}^2+\|K^{m,n}\ett_{[0,\eta_l]}\|_{\mcS^2_\bbQ}^2\leq C
\end{align}
for all $m,n\geq 0$ and $\bbQ\in\PrM^{L^z}$.
\end{lem}

\noindent\emph{Proof.} We note that $f^{m,n}$ satisfies the conditions of Proposition 3.1 in~\cite{HamLepIH}, guaranteeing the existence of a unique triple $(Y^{m,n},Z^{m,n},K^{m,n})\in\mcS^{2}\times\mcH^2\times\mcS^2$ that solves \eqref{ekv:trunk-time}.


We define
\begin{align*}
\gamma_s:=\frac{f^{m,n}(s,Y^{m,n}_s,Z^{m,n}_s)-f^{m,n}(s,0,Z^{m,n}_s)}{Y^{m,n}_s}\ett_{[Y^{m,n}_s\neq 0]}
\end{align*}
and set $e_{t,s}:=e^{\int_t^s\gamma_rdr}$ and $e_t:=e_{0,t}$. By Assumption~\ref{ass:on-rbsde-IH} we have that $e^{-\int_0^{s\wedge\eta_{m\vee n}} \underline L^y_rdr }\leq e_{t,s}\leq e^{\int_0^{s\wedge\eta_{m\vee n}} \bar L^y_r dr}$. Applying Ito's formula to $e_tY^{m,n}_t$ gives that for any $\tau\in\mcT_{t}$ we have
\begin{align*}
e_tY^{m,n}_t&=e_{\tau} Y^{m,n}_\tau +\int_t^\tau e_{s}(f^{m,n}(s,Y^{m,n}_s,Z^{m,n}_s)-\gamma_sY^{m,n}_s)ds-\int_t^\tau e_{s}Z^{m,n}_sdW_s+\int_t^Te_{s}dK^{m,n}_s
\\
&=e_{\tau} Y^{m,n}_\tau +\int_t^\tau e_{s}(f^{m,n}(s,0,0)+\zeta_sZ^{m,n}_s)ds-\int_t^\tau e_{s}Z^{m,n}_sdW_s+\int_t^\tau e_{s}dK^{m,n}_s
\end{align*}
where
\begin{align*}
\zeta_s:=\frac{f^{m,n}(s,0,Z^{m,n}_s)-f^{m,n}(s,0,0)}{|Z^{m,n}_s|^2}(Z^{m,n}_s)^\top\ett_{[Z^{m,n}_s\neq 0]}.
\end{align*}
By the Girsanov theorem (see \eg Chapter 15 in \cite{CohenElliottBook}), and since $\zeta_s=0$, $\Prob$-a.s.~on $[m\vee n,\infty)$, it follows that the measure $\bbQ^{\zeta}$ defined as $d\bbQ^{\zeta}:=\mcE(\zeta*W)_{m\vee n}d\Prob$ is a probability measure under which the process $W^{\zeta}:=W-\int_0^\cdot\zeta_sds$ is a Brownian motion. Furthermore, we have
\begin{align}\label{ekv:Y-mn-with-stopping}
Y^{m,n}_t=e_{t,\tau} Y^{m,n}_\tau +\int_t^\tau e_{t,s}f^{m,n}(s,0,0)ds-\int_t^\tau e_{t,s}Z^{m,n}_sdW^{\zeta}_s+\int_t^\tau e_{t,s}dK^{m,n}_s.
\end{align}
By the Cauchy-Schwarz inequality we find that
\begin{align*}
  \E^{\bbQ^\zeta}\Big[\big(\int_0^{\infty}|Z^{m,n}_s|^2 ds\big)^{1/2}\Big]\leq \E\big[|\mcE(\zeta*W)_{m\vee n}|^2\big]^{1/2}\E\Big[\int_0^{\infty}|Z^{m,n}_s|^2 ds\Big]^{1/2}<\infty
\end{align*}
and we conclude that $\int_t^{\cdot}e_{t,s}Z^{m,n}_s dW_s$ is $\bbQ^{\zeta}$-uniformly integrable. This implies that taking conditional expectation while picking $\tau=\inf\{r\geq t:Y^{m,n}_r=S_r\}$ in \eqref{ekv:Y-mn-with-stopping} gives
\begin{align*}
Y^{m,n}_t&=\E^{\bbQ^{\zeta}}\Big[\ett_{[\tau<\infty]}e_{t,\tau}S_{\tau}+\int_t^{\tau} e_{t,s}f^{m,n}(s,0,0)ds\big|\mcF_t\Big]
\end{align*}
or
\begin{align}\nonumber
|Y^{m,n}_t|&\leq \E^{\bbQ^{\zeta}}\Big[\sup_{r\in[t,\infty)}e^{\int_t^r\bar L^y_s ds}(S_{r})^+ + \int_t^{m\vee n}e^{\int_t^s\bar L^y_r dr}|f(s,0,0)|ds\big|\mcF_t\Big]
\\
&\leq \bar K^S_t + \bar K^f_t\label{ekv:Y-mn-bound}
\end{align}
and the estimate for $Y^{m,n}$ follows as $\bar K^S + \bar K^f\in\mcK^{2,\bar L^y,L^z}$.\\

The second bound is obtained by making a slight alteration in the proof of Proposition 3.5 in~\cite{ElKaroui1}. We pick an arbitrary $\zeta\in\mcZ^{L^z}$ and set $\tilde f^{m,n}(t,y,z):=f^{m,n}(t,y,z)-\zeta_t z$. Then the triple $(Y^{m,n},Z^{m,n},K^{m,n})$ satisfies
\begin{align*}
  \begin{cases}
    Y^{m,n}_t=Y^{m,n}_{\eta_l}+\int_t^{\eta_l} \tilde f^{m,n}(s,Y^{m,n}_s,Z^{m,n}_s)ds-\int_t^{\eta_l} Z^{m,n}_s dW^\zeta_s+K^{m,n}_{\eta_l}-K^{m,n}_t,\\
    Y^{m,n}_t\geq S_t,\, \forall t\in [0,{\eta_{l}}] \quad\mbox{and}\quad \int_0^{\eta_{l}} \left(Y^{m,n}_t-S_t\right)dK^{m,n}_t=0,
  \end{cases}
\end{align*}
which is not a standard reflected BSDE of the type in Theorem~\ref{thm:ElKaroui-rbsde}, since the data is not necessarily adapted to the filtration generated by $W^\zeta$. Nevertheless, we can apply Ito's formula to $(Y^{m,n})^2$ and find that
\begin{align}\nonumber
|Y^{m,n}_0|^2+\int_0^{\eta_l}|Z^{m,n}_s|^2ds &= |Y^{m,n}_{\eta_l}|^2+2\int_0^{\eta_l} Y^{m,n}_s \tilde f^{m,n}(s,Y^{m,n}_s,Z^{m,n}_s)ds
\\
&\quad-2\int_0^{\eta_l} Y^{m,n}_sZ^{m,n}_sdW^\zeta_s + 2\int_0^{\eta_l} Y^{m,n}_sdK^{m,n}_s\nonumber
\\
&\leq |Y^{m,n}_{\eta_l}|^2+\int_0^{\eta_l}(C|Y^{m,n}_s|^2+|f^{m,n}(s,0,0)|^2+\frac{1}{4}|Z^{m,n}_s|^2)ds\nonumber
\\
&\quad-2\int_0^{\eta_l}Y^{m,n}_sZ^{m,n}_sdW^{\zeta}_s + 2\int_0^{\eta_l} S_sdK^{m,n}_s,\label{ekv:ito-on-y2}
\end{align}
where we used the fact that $\tilde f^{m,n}$ is standard Lipschitz on $[0,\eta_l]$ (with coefficient $2l$) along with the relation $ab\leq \tfrac{1}{2}\kappa a^2+\tfrac{1}{2\kappa}b^2$ for a suitable $\kappa>0$, and the identity $\int_0^{\eta_{l}} \left(Y^{m,n}_t-S_t\right)dK^{m,n}_t=0$ to arrive at the inequality. The Burkholder-Davis-Gundy inequality and the relation $ab\leq \tfrac{1}{2}\kappa a^2+\tfrac{1}{2\kappa}b^2$ gives that
\begin{align*}
  \E^{\bbQ^{\zeta}}\Big[-2\int_0^{\eta_l}Y^{m,n}_sZ^{m,n}_sdW^{\zeta}_s\Big]&\leq C\E^{\bbQ^{\zeta}}\Big[(\int_0^{\eta_l}|Y^{m,n}_sZ^{m,n}_s|^2ds)^{1/2}\Big]
  \\
  &\leq \E^{\bbQ^{\zeta}}\Big[C\sup_{s\in[0,\eta_l]}|Y^{m,n}_s|^2+\frac{1}{4}\int_0^{\eta_l}|Z^{m,n}_s|^2ds\Big].
\end{align*}
Combining this with \eqref{ekv:ito-on-y2} we find that
\begin{align*}
  \E^{\bbQ^{\zeta}}\Big[\int_0^{\eta_l}|Z^{m,n}_s|^2ds\Big]&\leq C\E^{\bbQ^{\zeta}}\Big[\sup_{s\in[0,\eta_l]}|Y^{m,n}_s|^2+\int_0^{\eta_l}|f^{m,n}(s,0,0)|^2ds + \int_0^{\eta_l} S_sdK^{m,n}_s\Big]
  \\
  &\leq C(1+\E^{\bbQ^{\zeta}}\Big[\int_0^{\eta_l}|f^{m,n}(s,0,0)|^2ds + \int_0^{\eta_l} S_sdK^{m,n}_s\Big]),
\end{align*}
where the constant $C>0$ does not depend on $m,n$. By repeating the steps in the latter half of the proof of Proposition 3.5 in~\cite{ElKaroui1}, the bound in \eqref{ekv:unif-bound-Z-IH} now follows.\qed\\

\begin{lem}\label{lem:rbsde-trunk-stability}
For each $m,l\geq 0$ we have
\begin{align*}
\|(Y^{m,n'}-Y^{m,n})\ett_{[0,\eta_l]}\|_{\mcS^2}^2+\|(Z^{m,n'}-Z^{m,n})\ett_{[0,\eta_l]}\|_{\mcH^2}^2 + \|(K^{m,n'}-K^{m,n})\ett_{[0,\eta_l]}\|_{\mcS^2}^2\to 0,
\end{align*}
as $n,n'\to\infty$.
\end{lem}

\noindent\emph{Proof.} A trivial extension of the stability result, \eqref{ekv:ElKaroui-diff}, of Theorem~\ref{thm:ElKaroui-rbsde} to random terminal times implies that it suffice to show that $\E\big[|Y^{m,n'}_{\eta_l}-Y^{m,n}_{\eta_l}|^2\big]\to 0$ for all $l\geq 0$ as $n,n'\to\infty$. We assume, w.l.o.g., that $l\leq n\leq n'$ and have by \eqref{ekv:Y-mn-bound} that
\begin{align*}
|Y_{\eta_n}^{m,n}-Y_{\eta_n}^{m,n'}|\leq |Y_{\eta_n}^{m,n}|+|Y_{\eta_n}^{m,n'}|\leq \bar K^f_{\eta_n}+\bar K^S_{\eta_n}.
\end{align*}
For ease of notation we omit the superscripts $m,n$ and $m,n'$ and write $(Y',Z',K')$ for $(Y^{m,n'},Z^{m,n'},K^{m,n'})$ and $f'$ for $f^{m,n'}$. For $t\in[0,\eta_n]$, we have
\begin{align*}
Y'_{t}-Y_{t} &= Y'_{\eta_n}-Y_{\eta_n}+\int_t^{\eta_n} (f'(s,Y'_s,Z'_s)-f(s,Y_s,Z_s))ds
\\
&\quad-\int_t^{\eta_n} (Z'_s-Z_s)dW_s+K'_{\eta_n}-K'_t-(K_{\eta_n}-K_t).
\end{align*}
Now, let
\begin{align*}
\gamma_s:=\frac{f'(s,Y'_s,Z'_s)-f'(s,Y_s,Z'_s)}{Y'_s-Y_s}\ett_{[Y'_s\neq Y_s]}
\end{align*}
and
\begin{align*}
\zeta_s:=\frac{f'(s,Y_s,Z'_s)-f'(s,Y_s,Z_s)}{|Z'_s-Z_s|^2}(Z'_s-Z_s)^\top\ett_{[Z'_s\neq Z_s]}
\end{align*}
and set $e_{t,s}:=e^{\int_t^s\gamma_rdr}\leq e^{\int_0^s \bar L^y_r dr}<\infty$. Letting $\delta Y:=Y'-Y$, $\delta Z:=Z' - Z$ and $\delta f:=f'-f$ we get that for any $\tau\in\mcT_t$ we have
\begin{align*}
\delta Y_t&= e_{t,\tau}\delta Y_\tau+\int_t^\tau e_{t,s}(\delta f(s,Y_s,Z_s)+\zeta_s\delta Z_s)ds-\int_t^\tau e_{t,s}\delta Z_sdW_s+\int_t^\tau e_{t,s}d(\delta K)_s
\\
&= e_{t,\tau}\delta Y_\tau+\int_t^\tau e_{t,s}\delta f(s,Y_s,Z_s)ds-\int_t^\tau e_{t,s}\delta Z_sdW^{\zeta}_s+\int_t^\tau e_{t,s}d(\delta K)_s,
\end{align*}
where $W^{\zeta}=W-\int_0^\cdot\zeta_sds$. First, set $\tau=\inf\{r\geq t:Y_r=S_r\}\wedge \eta_n$ and we find that
\begin{align*}
\delta Y_t&=  e_{t,\tau}\delta Y_{\tau}+\int_t^{\tau} e_{t,s}\delta f(s,Y_s,Z_s)ds-\int_t^{\tau} e_{t,s}\delta Z_sdW^{\zeta}_s+\int_t^{\tau}e_{t,s}(dK'_s-dK_s)
\\
&\geq \ett_{[\tau=\eta_n]}e_{t,\eta_n}\delta Y_{\eta_n}+\int_t^{\tau} e_{t,s} \delta f(s,Y_s,Z_s)ds-\int_t^{\tau} e_{t,s}\delta Z_sdW^{\zeta}_s,
\end{align*}
where the inequality follows by noting that $\int_t^\tau dK_s=0$ and $\delta Y_\tau\ett_{[\tau<\eta_n]}=Y'_\tau-S_\tau\geq 0$. On the other hand, by picking $\tau=\inf\{r\geq t:Y'_r=S_r\}\wedge \eta_n$ we get
\begin{align*}
\delta Y_t&\leq \ett_{[\tau=\eta_n]}e_{t,\eta_n}\delta Y_{\eta_n}+\int_t^{\tau} e_{t,s}\delta f(s,Y_s,Z_s)ds-\int_t^{\tau} e_{t,s}\delta Z_sdW^{\zeta}_s.
\end{align*}
Now, on $[0,\eta_n]$ we have $f^{m,n}\equiv f^{m,n'}$ and taking the conditional expectation with respect to the measure $\bbQ^{\zeta}$, with $d\bbQ^{\zeta}=\mcE(\zeta*W)_n d\Prob$, we find that
\begin{align*}
|Y^{m,n'}_{\eta_l}-Y^{m,n}_{\eta_l}|&\leq \E^{\bbQ^{\zeta}}\big[e^{\int_0^{\eta_{n}} \bar L^y_r dr}|Y^{m,n'}_{\eta_n}-Y^{m,n}_{\eta_n}|\big|\mcF_{\eta_l}\big]
\\
&\leq \E^{\bbQ^{\zeta}}\big[e^{\int_0^{\eta_n}\bar L^y_r dr}(\bar K^f_{\eta_n}+\bar K^S_{\eta_n})\big|\mcF_{\eta_l}\big]
\end{align*}
In particular, this implies that
\begin{align*}
\E\big[|Y^{m,n'}_{\eta_l}-Y^{m,n}_{\eta_l}|^2\big]&\leq C\sup_{\bbQ\in\PrM^{L^z}}\E^\bbQ\big[|e^{\int_0^{\eta_n}\bar L^y_r dr}\bar K^f_{\eta_n}|^2+|e^{\int_0^{\eta_n}\bar L^y_r dr}\bar K^S_{\eta_n}|^2\big],
\end{align*}
where the right hand side tends to $0$ as $n\to\infty$.\qed\\

\begin{prop}\label{prop:rbsde-solu-IH}
The reflected BSDE \eqref{ekv:bsde-reflection-IH} admits a unique solution $(Y,Z,K)\in \mcS^{2,\bar L^y,{L^z}}\times \mcH^{2,loc}\times \mcS^{2,loc}$, with $K$ non-decreasing and $K_0=0$.
\end{prop}

\noindent\emph{Proof.} By comparison we find that for each $m\geq 0$ the sequence $(Y^{m,n})_{n\geq 0}$ of continuous processes is non-increasing and by Lemma~\ref{lem:rbsde-trunk-bound-IH} the sequence has bounded $\mcS^2$-norm. This implies that $Y^{m}:=\lim_{n\to\infty}Y^{m,n}$ exists, $\Prob$-a.s., and by Fatou's lemma the limit satisfies $\|e^{\int_0^\cdot \bar L^y_r dr}Y^{m}\|_{\mcS^2_\bbQ}^2\leq C$ for all $m\geq 0$ and $\bbQ\in\PrM^{L^z}$. Furthermore, we have by Lemma~\ref{lem:rbsde-trunk-stability} that
\begin{align*}
\|(Y^{m,n}-Y^{m,n'})\ett_{[0,\eta_l]}\|_{\mcS^2}\to 0
\end{align*}
for each $l\geq 0$, as $n,n'\to\infty$ implying that $Y^{m}$ is a continuous process.  Moreover, \eqref{ekv:Y-mn-bound} and Assumption~\ref{ass:on-rbsde-IH} implies that $Y^m\in\mcS^{2,\bar L^y,{L^z}}$. By \eqref{ekv:unif-bound-Z-IH} we also note that for each $l\geq 0$, $(\ett_{[0,\eta_l]}Z^{m,n})_{m,n\geq 0}$ is a uniformly bounded double-sequence in $\mcH^2$ and by again appealing to Lemma~\ref{lem:rbsde-trunk-stability} we have that
\begin{align*}
\|(Z^{m,n}-Z^{m,n'})\ett_{[0,\eta_l]}\|_{\mcH^2}\to 0,
\end{align*}
as $n,n'\to\infty$ and we conclude that there is a $Z^m$ and (similarly) also a $K^m$ such that $Z^{m,n}\ett_{[0,\eta_l]}\to Z^m\ett_{[0,\eta_l]}$ in $\mcH^p$ and $K^{m,n}\ett_{[0,\eta_l]}\to K^m\ett_{[0,\eta_l]}$, in $\mcS^2$-norm as $n\to\infty$. Now, by Theorem~\ref{thm:ElKaroui-rbsde} there is a unique triple $(\hat Y^m,\hat Z^m,\hat K^m)\in\mcS^2([0,\eta_l])\times\mcH^2([0,\eta_l])\times\mcS^2([0,\eta_l])$, with $\hat K^m$ non-increasing and $\hat K^m_0=0$, such that
\begin{align*}
  \begin{cases}
  \hat Y^m_t=Y^m_{\eta_l}+\int_t^{\eta_l} f^m(s,\hat Y^m_s,\hat Z^m_s)ds-\int_t^{\eta_l} \hat Z^m_sdW_s+\hat K^m_{\eta_l}-\hat K^m_t,\quad \forall t\in [0,{\eta_l}],\\
  \hat Y^m_t\geq S_t,\,\forall t\in [0,{\eta_l}]\quad\mbox{and}\quad \int_0^{\eta_l}(\hat Y^m_t-S_t)d\hat K^m_t=0,
  \end{cases}
\end{align*}
with $f^{m}:=\ett_{[0,\eta_m]}f^+-f^-$. By \eqref{ekv:ElKaroui-diff} it follows that
\begin{align*}
\|\hat Y^m-Y^{m,n}\|_{\mcS^2([0,\eta_l])}^2+\|\hat Z^m-Z^{m,n}\|_{\mcH^2([0,\eta_l])}^2 + \|\hat K^m-K^{m,n}\|_{\mcS^2([0,\eta_l])}^2\to 0,
\end{align*}
as $n\to\infty$ and by uniqueness of limits we conclude that $(Y^m_t,Z^m_t,K^m_t)_{0\leq t\leq \eta_n}=(\hat Y^m,\hat Z^m,\hat K^m)$ in $\mcS^2([0,\eta_l])\times\mcH^2([0,\eta_l])\times\mcS^2([0,\eta_l])$. Since for each $T\in[0,\infty)$ we have, outside of a $\Prob$-null set, that $\eta_l\geq T$ for $l\,(=l(\omega))$ sufficiently large, it follows that $(Y^m,Z^m,K^m)$ solves the reflected BSDE
\begin{align*}
  \begin{cases}
  Y^m_t=Y^m_T+\int_t^{T} f^m(s,Y^m_s,Z^m_s)ds-\int_t^{T} Z^m_sdW_s+K^m_{T}-K^m_t,\quad \forall t\in[0,T],\,\forall T\geq 0,\\
  Y^m_t\geq S_t,\,\forall t\in[0,\infty)\quad\mbox{and}\quad \int_0^{\infty}(Y^m_t-S_t)dK^m_t=0.
  \end{cases}
\end{align*}

If $(\tilde Y^m,\tilde Z^m,\tilde K^m)\in \mcS^{p,\bar L^y,{L^z}}\times \mcH^{2,loc}\times \mcS^{2,loc}$ is any other solution then for $l'\geq l$ we have
\begin{align*}
&\|(\tilde Y^m-Y^m)\ett_{[0,\eta_l]}\|_{\mcS^2}^2+\|(\tilde Z^m-Z^m)\ett_{[0,\eta_l]}\|_{\mcH^2}^2+\|(\tilde K^m-K^m)\ett_{[0,\eta_l]}\|_{\mcS^2}^2
\\
&\leq \E[\mcE(\gamma*W)_{\eta_l,\eta_{l'}}(|e^{\int_0^{\eta_{l'}}\bar L^y_r dr}\tilde Y^m_{\eta_{l'}}|^2+|e^{\int_0^{\eta_{l'}}\bar L^y_r dr}Y^m_{\eta_{l'}}|^2)],
\end{align*}
where $\zeta\in \mcZ^{L^z}$ and $C$ does not depend on $l'$. Since the right hand side tends to 0 as $l'\to\infty$, we conclude that $\tilde Y^m_t=Y^m_t$ and $\tilde K^m_t=K^m_t$, $\Prob$-a.s.~for all $t\in[0,\eta_l]$ and $\tilde Z^m\ett_{[0,\eta_l]}=Z^m\ett_{[0,\eta_l]}$, $d\Prob\times d\lambda$-a.e.~and uniqueness follows since $l\geq 0$ was arbitrary.\\

Now, arguing as in the proof of Lemma~\ref{lem:rbsde-trunk-stability} we find that for each $l\geq 0$, we have
\begin{align*}
\|(Y^{m'}-Y^{m})\ett_{[0,\eta_l]}\|_{\mcS^2}^2+\|(Z^{m'}-Z^{m})\ett_{[0,\eta_l]}\|_{\mcH^2}^2 + \|(K^{m'}-K^{m})\ett_{[0,\eta_l]}\|_{\mcS^2}^2\to 0,
\end{align*}
as $m,m'\to\infty$. By taking the limit as $m\to\infty$ and repeating the above arguments the desired result, thus, follows.\qed\\

\begin{rem}\label{rem:Z-in-mcH-loc-bbQ}
Having established that $Y\in \mcS^2_\bbQ$ for all $\bbQ\in\PrM^{L^z}$, we can use \eqref{ekv:ElKaroui-bound} (as we did in the proof of Lemma~\ref{lem:rbsde-trunk-bound-IH}) to conclude that for each $l\geq 0$, we have $Z\ett_{[0,\eta_l]}\in \mcH^2_\bbQ$ for all $\bbQ\in\PrM^{L^z}$.
\end{rem}

\bigskip

\begin{cor}\label{cor:rbsde-char-IH}
If $(Y,Z,K)$ solves \eqref{ekv:bsde-reflection-IH}, then $Y$ can be interpreted as the Snell envelope in the following way
\begin{equation}\label{ekv:Y-is-Snell-IH}
      Y_t=\esssup_{\tau\in\mcT_t}\E\bigg[\int_t^{\tau\wedge\eta_l} f(s,Y_s,Z_s)ds+\ett_{[\tau<\eta_l]}S_\tau+\ett_{[\tau\geq \eta_l]}Y_{\eta_l}\Big|\mcF_t\bigg]
\end{equation}
for all $l\geq 0$. In particular, with $D_t:=\inf\{s\geq t: Y_s=S_s\}$ we have the representation
\begin{equation}\label{ekv:Snell-attained}
      Y_t=\E\bigg[\int_t^{D_t\wedge\eta_l} f(s,Y_s,Z_s)ds+\ett_{[D_t<\eta_l]}S_{D_t}+\ett_{[D_t\geq\eta_l]}Y_{\eta_l}\Big|\mcF_t\bigg]
\end{equation}
for all $l\geq 0$ and $K_{D_t}-K_t=0$, $\Prob$-a.s.
\end{cor}

\noindent\emph{Proof.} The first two statements, \ie \eqref{ekv:Y-is-Snell-IH} and \eqref{ekv:Snell-attained}, are immediate from Proposition~\ref{prop:rbsde-solu-IH} and Theorem~\ref{thm:ElKaroui-rbsde}. That $K_{D_t}-K_t=0$ follows by letting $l\to\infty$ in $K_{D_t\wedge\eta_l}-K_t=0$.\qed\\


\section{Application to robust optimal stopping\label{sec:robust-impulse}}
We now apply the main results of Section~\ref{sec:rbsde-IH} to find weakly optimal solutions to robust optimal stopping problems. In particular, we are interested in finding $(\tau^*,\alpha^*)\in\mcT\times\mcA$ such that\footnote{Since we are in a zero-sum setting we look for a saddle-point, $J(\tau,\alpha^*)\leq J(\tau^*,\alpha^*)\leq J(\tau^*,\alpha)$ and have that $J(\tau^*,\alpha^*)=\sup_{\tau\in\mcT}\inf_{\alpha\in\mcA}J(\tau,\alpha)=\inf_{\alpha\in\mcA}\sup_{\tau\in\mcT}J(\tau,\alpha)$.}
\begin{align}
J(\tau^*,\alpha^*)=\sup_{\tau\in\mcT}\inf_{\alpha\in\mcA}J(\tau,\alpha).\label{ekv:maxmin}
\end{align}

Throughout, we assume the following forms on the drift and volatility terms in the forward SDE \eqref{ekv:forward-sde},
\begin{align*}
a(t,x,\alpha)=\left[\begin{array}{c}a_1(t,x)\\ a_2(t,x,\alpha)\end{array}\right] \quad{\rm and}\quad \sigma(t,x)=\left[\begin{array}{cc}\sigma_{1,1}(t,x)& 0\\ \sigma_{2,1}(t,x) & \sigma_{2,2}(t,x)\end{array}\right],
\end{align*}
where $a_1$ and $\sigma$ are functional Lipschitz in $x$ and $\sigma_{2,2}:[0,\infty)\times \bbD\to\R^{m_2\times m_2}$ (we let $\bbD$ denote the set of all \cadlag functions $x:[0,\infty)\to\R^d$ with the topology of uniform convergence) has an inverse, $\sigma_{2,2}^{-1}$, that is bounded on $[0,\infty)\times \bbD$.

For the purpose of solving \eqref{ekv:maxmin} we let $f$ be given by
\begin{align}
f(t,\omega,y,z):=\inf_{\alpha\in A}H(t,\omega,z,\alpha)=:H^{*}(t,\omega,z),\label{ekv:f-def}
\end{align}
where\footnote{We use the notation $(X_s)_{s\leq t}$ in arguments to emphasise that a function, for example, $\phi:[0,\infty)\times\bbD\times A\to \R$ at time $t$ only depend on the trajectory of $X$ on $[0,t]$.}
\begin{align*}
H(t,\omega,z,\alpha):=z\breve a(t,(X_s)_{s\leq t},\alpha)+e^{-\rho(t)}\phi(t,(X_s)_{s\leq t},\alpha),
\end{align*}
with
\begin{align*}
\breve a(t,x,\alpha):=\left[\begin{array}{c}0 \\ \sigma_{2,2}^{-1}(t,x)a_2(t,x,\alpha)\end{array}\right]
\end{align*}
and $X$ is the unique solution to the functional SDE
\begin{align}
X_t&=x_0+\int_0^t\tilde a(s,(X_r)_{ r\leq s})ds+\int_0^t\sigma(s,(X_r)_{ r\leq s})dW_s.\label{ekv:driftless-sde}
\end{align}
with
\begin{align*}
\tilde a(t,x):=\left[\begin{array}{c}a_1(t,x)\\ 0\end{array}\right].
\end{align*}

The approach we take to solve \eqref{ekv:maxmin} is to define a probability measure $\bbQ^{\alpha}$ under which $W^{\alpha}:=W-\int_0^\cdot\breve a(s,(X_r)_{ r\leq s},\alpha_s)ds$ is a Brownian motion and then let $\alpha=\alpha^*$, where $\alpha^*$ is a measurable selection of a minimizer to \eqref{ekv:f-def}. In particular, we note that for any $\alpha\in \mcA$, the 6-tuple $(\Omega,\mcF,\bbF,\bbQ^{\alpha},X,W^{\alpha})$ is a weak solution to \eqref{ekv:forward-sde} with control $\alpha$. In this regard our solution methodology renders a solution to \eqref{ekv:maxmin} in a weak sense.

However, before we move on to show optimality of this scheme, we give assumptions on $a,\sigma$ and $\phi,\psi$ under which the reflected BSDE \eqref{ekv:bsde-reflection-IH} with driver given by \eqref{ekv:f-def} attains a unique solution and $\bbQ^\alpha$ is a probability measure equivalent to $\Prob$.

\begin{ass}\label{ass:onSFDE}
For any $t,t'\geq 0$, $x,x'\in\bbD$ and $\alpha\in A$ and for some $q\geq 0$ we have:
\begin{enumerate}[i)]
  \item\label{ass:onSFDE-a-sigma} The coefficients $a:[0,\infty)\times\bbD\times A\to\R^{d}$ and $\sigma:[0,\infty)\times\bbD\to\R^{d\times d}$ are Borel-measurable, continuous in $t$ (and $\alpha$ when applicable) and satisfy the growth conditions
  \begin{align*}
    |a(t,(x_s)_{s\leq t},\alpha)|&\leq C^g_a(1+\sup_{s\leq t}|x_s|),
    \\
    |\sigma(t,(x_s)_{s\leq t})|&\leq C^g_\sigma(t),
  \end{align*}
  for some constant $C^g_a\geq 0$ and some continuous function $C^g_\sigma:[0,\infty)\to[0,\infty)$ of polynomial growth and the Lipschitz continuity
  \begin{align*}
    |a_1(t,(x_s)_{s\leq t})-a_1(t,(x'_s)_{s\leq t})|+|\sigma(t,(x_s)_{s\leq t})-\sigma(t,(x'_s)_{s\leq t})|&\leq C\sup_{s\leq t}|x'_s-x_s|.
  \end{align*}
  Moreover, $\sigma_{22}$ has a uniformly bounded matrix-inverse $\sigma_{22}^{-1}$.
  \item\label{ass:onSFDE-phi} The running reward $\phi:[0,\infty)\times \bbD\times A\to\R$ is $\mcB([0,\infty)\times\bbD\times A)$-measurable, continuous in $\alpha$ and satisfies the growth condition
  \begin{align*}
    |\phi(t,(x_s)_{s\leq t},\alpha)|\leq C^g_\phi(1+\sup_{s\leq t}|x_s|^q)
  \end{align*}
  for some $C^g_\phi>0$.
  \item\label{ass:onSFDE-psi} The terminal reward $\psi:[0,\infty)\times \bbD\to\R$ is $\mcB([0,\infty)\times\bbD)$-measurable with $t\mapsto \psi(t,(x_s)_{s\leq t})$ continuous for all $x\in\bbD$ and satisfies the growth condition
  \begin{align*}
    |\psi(t,(x_s)_{s\leq t})|\leq C^g_\psi(1+\sup_{s\leq t}|x_s|^q)
  \end{align*}
  for some $C^g_\psi>0$.
\end{enumerate}
\end{ass}

\bigskip

The $z$-coefficient of $H$ is a process $L^\alpha$ given by $L^\alpha_t:=\breve a(t,(X_s)_{s\leq t},\alpha)$ for all $\alpha\in A$. We thus note that $f$ as defined above has stochastic Lipschitz coefficient $L$, with $L_t:=\sup_{\alpha\in A}|L^a_t|\leq k_L(1+\sup_{s\in[0,t]}|X_s|)$ for some constant $k_L>0$. Together with the above assumptions this gives us the following important result:
\begin{lem}\label{lem:can-Girsanov}
Suppose that $\zeta\in \mcZ^L$, then $(\mcE(\zeta*W)_t)_{t\geq 0}$ is a martingale and for each $T\geq 0$ the measure $\bbQ^\zeta$ defined by $d\bbQ^\zeta=\mcE(\zeta*W)_Td\Prob$ is a probability measure under which $W^\zeta:=W-\int_0^\cdot\ett_{[0,T]}(s)\zeta_s ds$ is a Brownian motion.
\end{lem}

\noindent\emph{Proof.} See for example Lemma 0, pp. 456 in~\cite{Benes1971} or \cite{KarShreve1}, pp. 191.\qed\\

\subsection{Some preliminary estimates}
For any $\mcP_\bbF$-measurable process $\zeta$ we let $(\tilde X^\zeta,\tilde W)$ solve
\begin{align}
\tilde X^{\zeta}_t=x_0+\int_0^t(\tilde a(s,(\tilde X^{\zeta}_r)_{ r\leq s})+\sigma(s,(\tilde X^{\zeta}_r)_{ r\leq s})\zeta_s)ds + \int_0^t\sigma(s,(\tilde X^{\zeta}_r)_{ r\leq s})d\tilde W_s,\label{ekv:tilde-X-def}
\end{align}
whenever $\zeta$ is such that a solution exists. We will assume that $\tilde W$ is a $d$-dimensional $\bbF$-adapted Brownian motion under some auxiliary probability measure $\tilde \Prob$ on $(\Omega,\mcF)$, equivalent to $\Prob$. Using \eqref{ekv:tilde-X-def} we derive the following moment estimates:

\begin{prop}\label{prop:SFDEmoment}
Under Assumption~\ref{ass:onSFDE}, the SFDE \eqref{ekv:forward-sde} admits a weak solution for each $\alpha\in\mcA$. Furthermore, the solution has moments of all orders on compacts, in particular we have for $p\geq 2$, that
\begin{align}\label{ekv:SFDEmoment}
\sup_{\bbQ\in\PrM^L}\E^{\bbQ}\big[\sup_{t\in[0,T]}|X_t|^{p}\big]\leq C(1+e^{C^{g,p}_X(T)})
\end{align}
and
\begin{align}\label{ekv:SFDEmoment2}
\esssup_{\bbQ\in\PrM^L}\E^{\bbQ}\big[\sup_{s\in[t,T]}|X_s|^{p}\big|\mcF_t\big]\leq Ce^{C^{g,p}_X(T)-C^{g,p}_X(t)}(1+\sup_{s\in[0,t]}|X_s|^{p})
\end{align}
for all $0\leq t\leq T$, whenever $C^{g,p}_X(t)=(pC_a^g+\epsilon)t$ for some $\epsilon>0$.
\end{prop}

\noindent\emph{Proof.} Existence and uniqueness of solutions to \eqref{ekv:driftless-sde} is standard (see~\eg Chapter V of~\cite{Protter}). From Lemma~\ref{lem:can-Girsanov} it now follows that \eqref{ekv:forward-sde}, or more generally \eqref{ekv:tilde-X-def}, admits a weak solution whenever $\zeta$ is $\mcP_\bbF$-measurable and satisfies
\begin{align}\label{ekv:zeta-bound}
|\tilde a(s,(\tilde X^{\zeta}_r)_{ r\leq s})+\sigma(s,(\tilde X^{\zeta}_r)_{ r\leq s})\zeta_s|\leq C^g_a(1+\sup_{r\in[0,s]}|\tilde X^{\zeta}_s|)
\end{align}
for all $s\in[0,\infty)$.\\

To get the moment estimates we note that under $\bbQ^\zeta$ the dynamics of $X$ follow those of $\tilde X^\zeta$ with $\tilde W$ a Brownian motion. For $0\leq t\leq u$ we have, whenever $\zeta\equiv 0$ on $[0,t)$ and \eqref{ekv:zeta-bound} holds for all $s\in[t,\infty)$, that
\begin{align*}
|\tilde X^{\zeta}_u|\leq |\tilde X^{0}_t|+\int_t^uG^g_a(1+\sup_{r\in[0,s]}|\tilde X^\zeta_r|)ds + \Big|\int_t^u\sigma(s,(\tilde X^{\zeta}_r)_{ r\leq s})d\tilde W_s\Big|.
\end{align*}
By Gr\"onwall's inequality we get that
\begin{align}\label{ekv:tildeX-bound}
\sup_{r\in[t,u]}|\tilde X^{\zeta}_r|\leq Ce^{C^g_s(u-t)}(u-t+\sup_{r\in[0,t]}|\tilde X^{0}_r| +\sup_{r\in[t,u]}\Big|\int_t^r\sigma(s,(\tilde X^{\zeta}_r)_{ r\leq s})d\tilde W_s\Big|).
\end{align}
Both estimates now follow by the Burkholder-Davis-Gundy inequality (BDG for short) with $C^{g,q}_X(t)=qG^g_at+\epsilon t$ for any $\epsilon>0$.\qed\\

\subsection{The corresponding reflected BSDE\label{sec:robust-rbsde-FH}}
In the present section we show that there is a unique triple $(Y,Z,K)\in \mcS^{2,0,L}\times \mcH^{2,loc}\times \mcS^{2,loc}$ that solves the reflected BSDE
\begin{align}\label{ekv:bsde-robust-IH}
\begin{cases}
  Y_t=Y_T+\int_t^T H^{*}(s,Z_s)ds-\int_t^T Z_sdW_s+ K_T-K_t,\quad\forall t\in[0,T],\,\forall T>0, \\
  Y_T\to 0,\, \Prob\mbox{-a.s.~as }T\to\infty,\\
  Y_t\geq e^{-\rho(t)}\psi(t,(X_s)_{ s\leq t}),\,\forall t\in[0,\infty) \quad\mbox{and}\quad \int_0^\infty(Y_t-e^{-\rho(t)}\psi(t,(X_s)_{ s\leq t}))dK_t=0.
\end{cases}
\end{align}
Then, we will leverage the result in Corollary~\ref{cor:rbsde-char-IH} to find a weak solution to the robust optimal stopping problem in infinite horizon.

\begin{ass}\label{ass:on-rho}
The function $\rho:[0,\infty)\to[0,\infty)$, with $\rho(0)=0$, is continuously differentiable and such that $\tilde \rho(t):=qC^g_a t-\rho(t)$ satisfies $\tilde \rho(t)-\tilde \rho(s)\leq -\epsilon (t-s)$, for some $\epsilon>0$ and all $0\leq s\leq t$.
\end{ass}

We let $(\bar K^f_t)_{t\geq 0}$ and $(\bar K^S_t)_{t\geq 0}$ be continuous versions of
\begin{align}\label{ekv:Kf-def}
\bar K^f_t:=C^g_\phi\esssup_{\bbQ\in\PrM^L}\E^{\bbQ}\Big[\int_t^\infty e^{-\rho(s)}(1+\sup_{r\in[0,s]}|X_r|^{q})ds\Big|\mcF_t\Big]
\end{align}
and
\begin{align}\label{ekv:KS-def}
\bar K^S_t:=C^g_\psi\esssup_{\bbQ\in\PrM^L}\E^{\bbQ^\zeta}\Big[\sup_{s\in[t,\infty)}e^{-\rho(s)}(1+\sup_{r\in[0,s]}|X_r|^{q})\Big|\mcF_t\Big].
\end{align}
Note that the right hand sides of \eqref{ekv:Kf-def} and \eqref{ekv:KS-def} are both $\Prob$-supermartingales and thus have left and right limits, along rationals, everywhere. Arguing as in the proof of Lemma 4 in~\cite{Zamfirescu06} the existence of continuous versions follows.

Clearly, we have that $\bar K^f$ and $\bar K^S$ provide bounds for $f$ and $S$ in the manner described in Assumption~\ref{ass:on-rbsde-IH}. Moreover, combining Assumption~\ref{ass:onSFDE} and Assumption~\ref{ass:on-rho} gives us the following:

\begin{lem}\label{lem:K-f-K-S-valid}
The processes $\bar K^f$ and $\bar K^S$ both belong to $\mcK^{2,0,L}$. Moreover, there are constants $K^f>0$ and $K^S>0$ such that $\|\bar K^f\|^2_{\mcS^2_{\bbQ}}\leq K^f$ and $\|\bar K^S\|^2_{\mcS^2_{\bbQ}}\leq K^S$ for all $\bbQ\in\PrM^L$.
\end{lem}

\noindent\emph{Proof.} For $q=0$ the statement is immediate and so we only need to consider the case $q>0$. By \eqref{ekv:tildeX-bound} we find that whenever $\zeta_s=0$ on $[0,t)$ and \eqref{ekv:zeta-bound} holds for $s\geq t$ then for $s\geq t$ we have
\begin{align*}
e^{-\rho(s)}\sup_{r\in[t,s]}|\tilde X^{\zeta}_r|^{q}\leq Ce^{qC^g_s(s-t)-\rho(s)}((s-t)^q+\sup_{r\in[0,t]}|\tilde X^{0}_r|^q +\sup_{r\in[t,s]}\Big|\int_t^r\sigma(u,(\tilde X^{\zeta}_r)_{ r\leq u})d\tilde W_u\Big|^q)
\end{align*}
Using Fatou's lemma, Fubini's theorem and the BDG inequality we thus find that
\begin{align*}
\bar K^f_t&\leq Ce^{-\rho(t)}\lim_{T\to\infty}\E\Big[\int_t^Te^{-\epsilon(s-t)}(1+(s-t)^q+\sup_{r\in[0,t]}|\tilde X^{0}_r|^q +\sup_{r\in[t,s]}\Big|\int_t^r\sigma(s,(\tilde X^{\zeta}_r)_{ r\leq s})d\tilde W_s\Big|^q)ds\Big|\mcF_t\Big]
\\
&\leq Ce^{-\rho(t)}(1+\sup_{r\in[0,t]}|\tilde X^{0}_r|^q+\lim_{T\to\infty}\int_t^Te^{-\epsilon(s-t)}\Big(\int_t^s|C^g_\sigma(r)|^2dr\Big)^{q/2}ds)
\\
&\leq Ce^{-\rho(t)}(\sup_{r\in[0,t]}|\tilde X^{0}_r|^q+g(t)),
\end{align*}
where $g(t):=1+\int_t^\infty e^{-\epsilon(s-t)}\big(\int_t^s|C^g_\sigma(r)|^2dr\big)^{q/2}ds$ is of at most polynomial growth. In particular, since $-\rho(t)$ is non-increasing and bounded from above by $-\epsilon t$, we find that
\begin{align}
\sup_{s\in[t,\infty)}\bar K^f_s&\leq C(e^{-\frac{\epsilon}{2} t}+e^{-\rho(t)}\sup_{r\in[0,\infty)}e^{\rho(t)-\rho(r\vee t)}|X_r|^q).\label{ekv:Kf-sup-bound}
\end{align}
To show that $\bar K^f\in\mcS^2_\bbQ$ for all $\bbQ\in \PrM^L$ we need to show that the right hand side is square $\E^\bbQ$-integrable. We thus go back and assume that $\zeta$ is such that \eqref{ekv:zeta-bound} holds for all $s\in[0,\infty)$. With $h(r):=(\rho(t)-\rho(r\vee t))/q$ we then have for $0\leq t\leq s$ that
\begin{align*}
e^{h(s)}\tilde X^{\zeta}_s=\tilde X^{\zeta}_t+\int_t^s e^{h(r)}(\tilde a(r,(\tilde X^{\zeta}_u)_{u\leq r})+\sigma(r,(\tilde X^{\zeta}_u)_{ u\leq r})\zeta_r+h'(r)\tilde X^{\zeta}_r)dr + \int_t^s e^{h(r)}\sigma(r,(\tilde X^{\zeta}_u)_{ u\leq r})d\tilde W_r
\end{align*}
and since $h'(r)< -C^g_a$ we find that
\begin{align}\label{ekv:e-X-bound}
\sup_{r\in[0,s]}e^{h(r)}|\tilde X^{\zeta}_r|\leq \sup_{r\in[0,t]}|\tilde X^{\zeta}_r|+\sup_{t'\in[t,s]}\Big|\int_t^{t'}e^{h(r)}\sigma(r,(\tilde X^{\zeta}_u)_{ u\leq r})d\tilde W_r\Big|.
\end{align}
In particular, setting $t=0$ the BDG inequality now gives that for all $\bbQ\in\PrM^L$, we have
\begin{align*}
\|\bar K^f\|_{\mcS^2_\bbQ}^2\leq C(1+\big(\int_0^\infty |e^{h(s)}C^g_\sigma(s)|^{2} ds\big)^q)=:K_f<\infty.
\end{align*}
Moreover, for any $T\geq 0$ and $\bbQ\in\PrM^L$ we find by combining \eqref{ekv:Kf-sup-bound} and \eqref{ekv:e-X-bound} with \eqref{ekv:SFDEmoment} that
\begin{align*}
\|\ett_{[T,\infty)}\bar K^f\|_{\mcS^2_\bbQ}^2\leq Ce^{- \frac{\epsilon}{2} T}(1+\big(\int_0^\infty |e^{h(s)}C^g_\sigma(s)|^{2} ds\big)^q),
\end{align*}
where $C>0$ does not depend on $T$ and $\bbQ$ and it follows that $\bar K^f\in\mcK^{2,0,L}$. The result for $\bar K^S$ is immediate from the above.\qed\\

This leads us to the following conclusion:
\begin{prop}\label{prop:robust-rbsde-has-solution}
Under Assumption \ref{ass:onSFDE} and Assumption \ref{ass:on-rho}, the reflected BSDE \eqref{ekv:bsde-robust-IH} admits a unique solution $(Y,Z,K)\in \mcS^{2,0,L}\times \mcH^{2,loc}\times \mcS^{2,loc}$.
\end{prop}

\noindent\emph{Proof.} By Lemma~\ref{lem:K-f-K-S-valid} we have that the driver $f(t,y,z):=H^{*}(t,z)$ and the obstacle $S_t=e^{-\rho(t)}\psi(t,(X_s)_{ s\leq t})$ satisfy the conditions of Assumption~\ref{ass:on-rbsde-IH} by which Proposition~\ref{prop:rbsde-solu-IH} implies the existence of a unique solution to \eqref{ekv:bsde-robust-IH}.\qed\\


\subsection{Robust optimal stopping in infinite horizon}

By Benes’ selection Theorem (\cite{Benes1971}, Lemma 5, pp. 460), there exists a $\mcP\otimes\mcB(\R^d)/\mcB(A)$-measurable function $\alpha(t,\omega,z)$ such that for any given $(t,\omega, z) \in [0,\infty)\times\Omega\times \R^{d}$, we have
\begin{align*}
H(t,\omega,z,\alpha(t,\omega,z))=\inf_{\alpha\in A}H(t,\omega,z,\alpha),
\end{align*}
$\Prob$-a.s.

The following theorem shows that we can extract the optimal pair $(\tau^*,\alpha^*)$ from the map $(t,\omega,z)\mapsto \alpha(t,\omega,z)$ and the solution to \eqref{ekv:bsde-robust-IH}.

\begin{thm}
Under Assumption \ref{ass:onSFDE} and Assumption \ref{ass:on-rho}, let the triple $(Y,Z,K)\in \mcS^{2,0,L}\times \mcH^{2,loc}\times \mcS^{2,loc}$ solve \eqref{ekv:bsde-robust-IH}. Then the pair $(\tau^*,\alpha^*)\in\mcT\times \mcA$, defined as:
\begin{align*}
\tau^*:=\inf \{t \geq 0:\:Y_t=e^{-\rho(t)}\psi(t,(X_s)_{s\leq t})\}
\end{align*}
and
\begin{align*}
\alpha^*_t:=\alpha(t,Z_t)
\end{align*}
is an optimal pair in the sense that
\begin{align}\label{ekv:is-robust-IH}
Y_0=J(\tau^*,\alpha^*)=\sup_{\tau\in\mcT}\inf_{\alpha\in\mcA}J(\tau,\alpha).
\end{align}
\end{thm}

\noindent\emph{Proof.} Since $Y_0$ is $\mcF_0$-measurable and $\mcF_0$ is trivial and $Z\ett_{[0,\eta_l]}\in\mcH^{2}_{\bbQ}$ for all $\bbQ\in\PrM^L$ and all $l\geq 0$ (see Remark~\ref{rem:Z-in-mcH-loc-bbQ}), we have by Corollary~\ref{cor:rbsde-char-IH} that
\begin{align*}
Y_0&=\E^{\bbQ^{\alpha^*}_l}\Big[\ett_{[\tau^*< {\eta_l}]}e^{-\rho(\tau^*)}\psi(\tau^*,(X_s)_{ s\leq \tau^*})+\ett_{[\tau^*\geq {\eta_l}]}Y_{\eta_l}+\int_0^{\tau^*\wedge\eta_l} H^{*}(s,Z_s)ds
-\int_0^{\tau^*\wedge\eta_l} Z_sdW_s\Big]
\\
&=\E^{\bbQ^{\alpha^*}_l}\Big[\ett_{[\tau^*< {\eta_l}]}e^{-\rho(\tau^*)}\psi(\tau^*,(X_s)_{ s\leq \tau^*})+\ett_{[\tau^*\geq {\eta_l}]}Y_{\eta_l}+\int_0^{\tau^*\wedge\eta_l}  e^{-\rho(s)}\phi(s,(X_r)_{ r\leq s},\alpha^*_s)ds\Big],
\end{align*}
where $\bbQ^{\alpha}_l$ is the measure, equivalent to $\Prob$, under which $W^{\alpha}:=W-\int_0^{\cdot\wedge\eta_l} L^{\alpha}_sds$ is a martingale. Now, by Proposition~\ref{prop:robust-rbsde-has-solution}, $Y\in\mcS^{2,0,L}$ and for any $(\tau,\alpha)\in\mcT\times\mcA$ we have
\begin{align*}
\sup_{\bbQ\in\PrM^L}\E^{\bbQ}\Big[\ett_{[\tau\geq {\eta_l}]}e^{-\rho(\tau)}|\psi(\tau,(X_s)_{ s\leq {\tau}})|+\int_{\eta_l}^{\infty}e^{-\rho(s)}|\phi(s,(X_r)_{ r\leq s},\alpha_s)|ds\Big]\to 0,
\end{align*}
as $l\to\infty$ which leads us to conclude that
\begin{align*}
Y_0=J(\tau^*,\alpha^*).
\end{align*}
Moreover, for any other $\tau\in\mcT$ we note that for all $l\geq 0$ we have
\begin{align*}
Y_0&\geq\E^{\bbQ^{\alpha^*}_l}\Big[\ett_{[\tau< {\eta_l}]}e^{-\rho(\tau)}\psi(\tau,(X_s)_{ s\leq \tau})+\ett_{[\tau\geq {\eta_l}]}Y_{\eta_l}+\int_0^{\tau\wedge\eta_l} H^{*}(s,Z_s)ds
-\int_0^{\tau\wedge\eta_l} Z_sdW_s\Big]
\end{align*}
where the right hand side converges, decreasingly, to $J(\tau,\alpha^*)$ as $l\to\infty$ and we conclude that
\begin{align*}
Y_0=J(\tau^*,\alpha^*)=\sup_{\tau\in\mcT}J(\tau,\alpha^*).
\end{align*}
To realize that $\alpha^*$ is an optimal response we note that for any $\alpha \in\mcA$, we have
\begin{align*}
Y_0&=\E^{\bbQ^{\alpha}_l}\Big[\ett_{[\tau^{*}< {\eta_l}]}e^{-\rho(\tau^{*})}\psi(\tau^{*},(X_s)_{ s\leq \tau^{*}})+\ett_{[\tau^{*}\geq {\eta_l}]}Y_{\eta_l}+\int_0^{\tau^{*}\wedge\eta_l} H^{*}(s,Z_s)ds
-\int_0^{\tau^{*}\wedge\eta_l} Z_sdW_s\Big]
\\
&=\E^{\bbQ^{\alpha}_l}\Big[\ett_{[\tau^{*}< {\eta_l}]}e^{-\rho(\tau^{*})}\psi(\tau^{*},(X_s)_{ s\leq \tau^{*}})+\ett_{[\tau^{*}\geq {\eta_l}]}Y_{\eta_l}+\int_0^{\tau^{*}\wedge\eta_l}  e^{-\rho(s)}\phi(s,(X_r)_{ r\leq s},\alpha_s)ds\Big]
\\
&\quad+\E^{\bbQ^{\alpha}_l}\Big[\int_0^{\tau^{*}\wedge\eta_l}(H^{*}(s,Z_s)-H(s,Z_s,\alpha_s))ds\Big]
\\
&\leq \E^{\bbQ^{\alpha}_l}\Big[\ett_{[\tau^{*}< {\eta_l}]}e^{-\rho(\tau^{*})}\psi(\tau^{*},(X_s)_{ s\leq \tau^{*}})+\ett_{[\tau^{*}\geq {\eta_l}]}Y_{\eta_l}+\int_0^{\tau^{*}\wedge\eta_l}  e^{-\rho(s)}\phi(s,(X_r)_{ r\leq s},\alpha_s)ds\Big].
\end{align*}
Repeating the above argument we find that the right-hand side tends to $J(\tau^{*},\alpha)$ as $l\to\infty$ and we conclude that $J(\tau^{*},\alpha^*)\leq J(\tau^{*},\alpha)$. In particular, this proves that $(\tau^*,\alpha^*)\in\mcT\times \mcA$ is a saddle point in the sense that $J(\tau,\alpha^*)\leq J(\tau^*,\alpha^*)\leq J(\tau^*,\alpha)$ for any $(\tau^*,\alpha^*)\in\mcT\times \mcA$. The second equality in \eqref{ekv:is-robust-IH} is an immediate consequence of this.\qed\\

\bibliographystyle{plain}
\bibliography{RBSDEih_ref}
\end{document}